\documentclass{amsart}

\newtheorem*{theorem}{Theorem}
\theoremstyle{definition}
\newtheorem*{definition}{Definition}

\begin{document}

\title[On the Pontryagin Hamiltonian for Autonomous Problems]{On the Constancy 
of the Pontryagin Hamiltonian for Autonomous Problems}

\author{Delfim F. M. Torres}
\address{Department of Mathematics \\
         University of Aveiro \\
         3810-193 Aveiro, Portugal}
\email{delfim@mat.ua.pt}

\date{March 3, 2003}


\begin{abstract}
We provide a new, simpler, and more direct proof of the well known fact
that for autonomous optimal control problems the Pontryagin extremals
evolve on a level surface of the respective Pontryagin Hamiltonian.
\end{abstract}

\maketitle


Given sets $\Omega \subseteq \mathbb{R}^r$
and $\mathcal{F} \subseteq \mathbb{R}^{2n}$,
constants $a < b$, and two continuous functions 
$L(x,u) : \mathbb{R}^n \times \Omega \rightarrow \mathbb{R}$ and
$\varphi(x,u) : \mathbb{R}^n \times \Omega \rightarrow
\mathbb{R}^n$ with continuous derivatives with respect 
to $x$, we define the autonomous optimal control problem 
as the minimization or maximization of the cost functional
$I\left[x(\cdot),u(\cdot)\right] = \int_a^b L\left(x(t),u(t)\right)
\mathrm{d}t$,
called the performance index,
among all the solutions of the vector differential equation
$\dot{x}(t) = \varphi\left(x(t),u(t)\right)$ for almost all
$t \in [a,b]$,
subject to the boundary conditions $\left(x(a),x(b)\right) \in \mathcal{F}$.
The \emph{state trajectory} $x(\cdot)$ is a $n$-vector absolutely continuous
function and the \emph{control} $u(\cdot)$ is a $r$-vector measurable and bounded
function satisfying the control constraint $u(t) \in \Omega$:
$x(\cdot) \in W_{1,1}\left([a,b];\mathbb{R}^n\right)$, 
$u(\cdot) \in L_{\infty}\left([a,b];\Omega\right)$.
The problem is denoted by $(P)$.
The celebrated Pontryagin maximum principle
\cite{MR29:3316b}, which is a first-order necessary optimality
condition for optimal control, provides a generalization
of the classical calculus of variations first-order necessary
optimality conditions. It asserts that the 
minimizers or maximizers of the optimal control problems
are to be found among the Pontryagin extremals.
\begin{definition}
Let us associate to the optimal control problem $(P)$ 
the Hamiltonian $H$ defined by $H(x,u,\psi_0,\psi) 
= \psi_0 L(x,u) + \psi \cdot \varphi(x,u)$.
A quadruple $\left(x(\cdot),u(\cdot),\psi_{0},\psi(\cdot)\right)$,
where $\psi_0 \le 0$ is a constant and $\psi(\cdot)$ a $n$-vector
absolutely continuous function with domain $[a,b]$,
is called a Pontryagin extremal if it
satisfies the control system
$\dot{x}(t) = \frac{\partial H}{\partial \psi}\left(x(t),u(t),\psi_0,\psi(t)\right)$;
the adjoint system $\dot{\psi}(t) = 
- \frac{\partial H}{\partial x}\left(x(t),u(t),\psi_0,\psi(t)\right)$;
and the maximality condition
$H\left(x(t),u(t),\psi_0,\psi(t)\right) = \max_{v \in \Omega}
H\left(x(t),v,\psi_0,\psi(t)\right)$.
\end{definition}
In the present note we are interested  in the following well-known 
result \cite{MR29:3316b}.
\begin{theorem}
If $\left(x(\cdot),u(\cdot),\psi_0,\psi(\cdot)\right)$ is
a Pontryagin extremal of $(P)$, then
\begin{equation}
\label{eq:HamilConst}
H(x(t),u(t),\psi_0,\psi(t)) \equiv constant \, , \quad t \in [a,b] \, .
\end{equation}
\end{theorem}
The Theorem has several important applications.
In classical mechanics, \eqref{eq:HamilConst} corresponds to conservation
of energy (\textrm{cf. e.g.} \cite{MR55:4815});
in economics to the Ramsey rule for optimal saving
or to the constancy of the welfare measure of national income
(\textrm{cf. e.g.} \cite{MR1118388});
while in the calculus of variations
it corresponds to the second Erdmann necessary optimality condition
(\textrm{cf. e.g.} \cite{MR91e:49001}).
Although the Theorem is a consequence of
the Pontryagin maximum principle, standard proofs use
sophisticated and lengthy extra arguments
(\textrm{cf. e.g.} \cite{MR51:8914,MR58:33350c,MR29:3316b}).
Here we introduce a new approach. We show
that equality \eqref{eq:HamilConst} can be trivially obtained
when one applies the Pontryagin maximum principle
to a suitable auxiliary optimal control problem. This new
problem is obtained introducing a new state variable, and
the conclusion follows from the fact that the new
Hamiltonian does not depend on the added state variable and therefore
the corresponding multiplier must be constant.
Our proof is motivated by Emmy Noether's theorem 
of optimal control \cite{delfimEJC}, and the use of
ignorable or kinosthenic variables in mechanics 
(\textrm{cf. e.g.} \cite{MR55:4815}).
\begin{proof}
It is clear that the autonomous problem of optimal
control is time-invariant: problem $(P)$ is invariant
under the parameter transformation $\tau = t + s$
(autonomous problems show symmetry under time translation). 
Our viewpoint is to consider the
parameter $s$ not as a constant but as a function of $\tau$.
We introduce a new state variable
$s(\tau) \in W_{1,\infty}\left([\alpha,\beta];\mathbb{R}\right)$
which satisfies the boundary conditions
$s(\alpha) = \alpha - a$, $s(\beta) = \beta - b$, and
whose derivative takes value on the open set $(0,1)$.
Doing the change of variable $t = \tau - s(\tau)$,
$\mathrm{d}t = \left(1 - s'(\tau)\right)\mathrm{d}\tau$,
and introducing the notation $z(\tau) = x\left(\tau-s(\tau)\right) = x(t)$,
$w(\tau) = u\left(\tau-s(\tau)\right) = u(t)$, 
we get from $(P)$ the following optimal control problem:
\begin{gather}
J\left[z(\cdot),w(\cdot),v(\cdot)\right] =
\int_{\alpha}^{\beta} L\left(z(\tau),w(\tau)\right) 
\left(1-v(\tau)\right) \mathrm{d}\tau
\longrightarrow \textrm{extr} \, , \notag \\
\begin{cases}
z'(\tau) = \varphi\left(z(\tau),w(\tau)\right) \left(1-v(\tau)\right) \, , \\
s'(\tau) = v(\tau) \, ,
\end{cases} \label{Prb:AuxP} \\
\left(z(\alpha),z(\beta)\right) \in \mathcal{F} \, , \quad
\left(s(\alpha),s(\beta)\right) = \left(\alpha - a,\beta - b\right) \, . \notag
\end{gather}
Compared to $(P)$, problem \eqref{Prb:AuxP} has one more
state variable and one more control variable. Namely,
the state variables are $z(\cdot)$ and $s(\cdot)$,
and the control variables are $w(\cdot)$ and $v(\cdot)$.
The Hamiltonian $\mathcal{H}$ associated to problem \eqref{Prb:AuxP} 
does not depend on the state variable $s$:
$\mathcal{H}\left(z,v,w,p_0,p_z,p_s\right)
= \left(p_0 L(z,w) + p_z \cdot \varphi(z,w)\right) (1-v) + p_s v
= H(z,w,p_0,p_z) (1-v) + p_s v$.
From the maximality condition it must be the case that
$\frac{\partial \mathcal{H}}{\partial v} = 0$, that is,
$p_s(\tau) = H\left(z(\tau),w(\tau),p_0,p_z(\tau)\right)$;
while from the adjoint system it follows that
$p_{s}'(\tau) = - \frac{\partial \mathcal{H}}{\partial s} = 0$.
One concludes that $H\left(z(\tau),w(\tau),p_0,p_z(\tau)\right) 
= \text{constant}$, and equality \eqref{eq:HamilConst} holds 
trivially from the fact that the set of Pontryagin extremals 
of \eqref{Prb:AuxP} is richer than that of $(P)$.
\end{proof}



\end{document}